\documentclass[10pt]{article}

\usepackage{amsmath}    % need for subequations
\usepackage{graphicx}   % need for figures
\usepackage{verbatim}   % useful for program listings
\usepackage{color}      % use if color is used in text
\usepackage{subfigure}  % use for side-by-side figures
\usepackage{amssymb,tikz}
\usepackage{soul}

\usepackage[bookmarks, bookmarksopen, bookmarksdepth=3, colorlinks=true,
linkcolor=black, anchorcolor=black, citecolor=black, filecolor=black,
menucolor=black, runcolor=black, urlcolor=black, pdfencoding=unicode
]{hyperref}

\newtheorem{theorem}{Theorem}[section]
\newtheorem{proposition}[theorem]{Proposition}

\newtheorem{corollary}[theorem]{Corollary}
\newtheorem{lemma}[theorem]{Lemma}

\newtheorem{problem}[theorem]{Problem}
\newtheorem{example}[theorem]{Example}

\newcommand{\proof}{\noindent{\bf Proof.\ }}
\newcommand{\qed}{\hfill $\square$\medskip}
\def\cp{\,\square\,}

\renewcommand{\gg}{\gamma_{g}}
\newcommand{\gtg}{\gamma_{tg}}
\newcommand{\gcg}{\gamma_{cg}}
\newcommand{\gc}{\gamma_{c}}
\newcommand{\sgcg}{\widetilde{\gamma}_{cg}}

\DeclareMathOperator {\diam} {diam}

\DeclareMathOperator {\CL} {CL}
\DeclareMathOperator {\ML} {ML}

\begin{document}
	
% Define style for nodes
\tikzstyle{every node}=[circle, draw, fill=black!10,
	inner sep=0pt, minimum width=4pt]

% Define style for nodes
\tikzstyle{every node}=[circle, draw, fill=black!10,
                        inner sep=0pt, minimum width=4pt]

\title{Connected domination game: predomination, Staller-start game, and
	lexicographic products}
 
 \author{
 	Vesna Ir\v si\v c $^{a, b}$
 }
 
 \date{\today}
 
 \maketitle
 % \vspace{-0.8 cm}
 \begin{center}
 	$^a$ Institute of Mathematics, Physics and Mechanics, Ljubljana, Slovenia\\
 	{\tt vesna.irsic@fmf.uni-lj.si }
 	\medskip
 
 	$^b$ Faculty of Mathematics and Physics, University of Ljubljana, Slovenia\\
 		
 \end{center}
 
 \begin{abstract}
	The connected domination game was recently introduced by Borowiecki, Fiedorowicz and Sidorowicz as another variation of the domination game. The rules are essentially the same, except that the set of played vertices must be connected at all stages of the game. We answer a problem from their paper regarding the relation between the number of moves in a game where Dominator/Staller starts the game. In this paper we also study the relation to the diameter and present graphs with small connected game domination number. We determine the values on the lexicographic product graphs, and consider the effect of predomination of a vertex on the connected game domination number.
%	The connected domination game was recently introduced by Borowiecki, Fiedorowicz and Sidorowicz as another variation of the domination game. We answer a problem from their paper regarding the relation between the number of moves in a game where Dominator/Staller starts the game. Additionally, we study the relation to the diameter and present graphs with small connected game domination number. We also determine the value on the lexicographic products, and consider the effect of predomination of a vertex.
 \end{abstract}
 
 
 \medskip
 \noindent{\bf AMS Subj.\ Class.:} 05C57, 05C69, 05C76

%%%%%%%%%%%%%%%%%%%%%%%%%%%
\section{Introduction}
\label{sec:intro}
%%%%%%%%%%%%%%%%%%%%%%%%%%%

The \emph{connected domination game} was introduced very recently in~\cite{connected} as a game played on a graph $G$ by two players, Dominator and Staller, with the following rules. Dominator's goal is to finish the game with the smallest possible number of moves, while Staller aims to prolong the game. The players alternately select vertices such that each move dominates at least one vertex that was not dominated by the previous moves, and the selected vertex is connected to at least one already played vertex. The game ends when no more moves are possible, i.e.\ the selected vertices form a connected dominating set.

If both players play optimally, the number of moves on a given graph is the \emph{connected game domination number} $\gcg(G)$ if Dominator starts the game on $G$. In this case we say for short that a \emph{D-game} is played on $G$. If Staller starts the game (shortly, \emph{S-game}), the number of moves is the 
\emph{Staller-start connected game domination number} $\gcg'(G)$.

The connected domination game is a variation of the classical domination game, which was introduced in 2010 by Bre\v{s}ar, Klav\v{z}ar and Rall~\cite{dom} and has been widely studied ever since (see for example~\cite{G-e, dom35, book, domEdgeCuts, domPathsCycles, domLargest}), along with the total domination game which was introduced in~\cite{totDom} (see also~\cite{3/4_2, totDomCP, 3/4_1, 3/4_3, trees}). Recently, three other variations of the game have been introduced~\cite{domGames}.

In the first paper on the connected domination game~\cite{connected}, the connected game domination number of trees and 2-paths is determined, along with an upper bound on 2-trees. The game is studied on the Cartesian product graphs. Some results about the Staller-start and Staller-pass game are also presented. In the Staller-pass game Staller is allowed to pass some moves instead of playing a vertex. The moves when she passes do not count into the number of moves in the game. The number of moves in a game on $G$ when Staller can pass at most $k$ times is $\hat{\gamma}_{cg}^k(G).$ It was proved in~\cite{connected} that $\gcg(G) \leq \hat{\gamma}_{cg}^k(G) \leq \gcg(G)+k$. 

Another important variation of the connected domination game presented in~\cite{connected} is the \emph{connected domination game with Chooser} which turns out to be an extremely useful tool for proving some bounds for the connected game domination number. The rules of the game are essentially the same as in the usual game, except that there is another player, Chooser, who can make zero, one or more moves after any move of Dominator or Staller. The rules for his move to be legal are the same as for Dominator and Staller. Note that Chooser has no specific goal, he can help either Dominator or Staller. For later use we also state the Chooser Lemma from~\cite{connected}.

\begin{lemma}[Chooser Lemma]
	\label{lema:chooser}
	Consider the connected domination game with Chooser on a graph $G$. Suppose that in the game Chooser picks $k$ vertices, and that both Dominator and Staller play optimally. Then at the end of the game the number of played vertices is at most $\gcg(G) + k$ and at least $\gcg(G) - k$.
\end{lemma} 

We now recall some basic definitions. The \emph{(closed) neighborhood} of a vertex $v$ in a graph $G$ is denoted by $N[v]$. The (closed) neighborhood of a set of vertices $S \subseteq V(G)$ is $N[S] = \bigcup_{v \in S} N[v]$. Similarly, the \emph{open neighborhood} of $v$ is $N(v) = N[v] \setminus \{v\}$. A set $S \subseteq V(G)$ is a \emph{dominating set} of the graph $G$ if $N[S] = V(G)$. The minimal cardinality of a dominating set is the \emph{domination number} $\gamma(G)$. The \emph{connected domination set} is a set $S$ that is both connected and dominating. The smallest cardinality of such set is the \emph{connected domination number} $\gc(G)$. Recall also that a \emph{join} of graphs $G$ and $H$ is a graph with vertex set $V(G) \cup V(H)$ and edges $E(G) \cup E(H) \cup \{ gh \; ; \; g \in G, h \in H \}$. Additionally, we denote $[n] = \{1, \ldots, n\}$, and $\Delta(G)$ as the maximum degree of vertices of $G$.

Let $G$ and $H$ be graphs. The \emph{lexicographic product} $G[H]$ of graphs $G$ and $H$ has the vertex set $V(G) \times V(H)$ and the edges $(g_1, h_1) \sim (g_2, h_2)$ if $g_1g_2 \in E(G)$, or $g_1=g_2$ and $h_1h_2 \in E(H)$. Note that $G[H]$ can be obtained by substituting each vertex $g$ of $G$ with a copy of $H$ (denoted by $H_g$) and connecting all vertices in $H_{g_1}$ with all vertices in $H_{g_2}$ if and only if $g_1g_2 \in E(G)$.

In this paper we first characterize graphs with small connected game domination number and observe a relation between the diameter of a graph and values of $\gcg(G)$, $\gcg'(G)$. Next we solve Problem 1 from~\cite{connected} regarding the difference between $\gcg(G)$ and $\gcg'(G)$. In Section~\ref{sec:lexicographic} we discuss the connected domination game on the lexicographic product of graphs. The last topic we study is the predomination of vertices and its effect on the connected game domination number. We conclude the paper with a discussion about possibilities how to define connected game domination critical graphs.

%Mogoče za omenit:\\
%- uni-cyclic graphs\\
%- bluffing?\\
%- 2-connected graphs?\\
%- Blok-grafi?

%%%%%%%%%%%%%%%%%%%%%%%%%%%%%%%%%%%%%%%%%%%%%%%%%%%%%%%%%%
%%%%%%%%%%%%%%%%%%%%%%%%%%%%%%%%%%%%%%%%%%%%%%%%%%%%%%%%%%
\section{Graphs with small connected game domination number and the relation to the diameter of a graph}
\label{sec:small+diam}
%%%%%%%%%%%%%%%%%%%%%%%%%%%%%%%%%%%%%%%%%%%%%%%%%%%%%%%%%%
%%%%%%%%%%%%%%%%%%%%%%%%%%%%%%%%%%%%%%%%%%%%%%%%%%%%%%%%%%

Motivated by~\cite{domSmall} where graphs with small game domination numbers were studied, we consider graphs with the connected domination number $1$~and~$2$. (The characterization of graphs with connected game domination number~$3$ seems rather technical, thus we omit it.)

\begin{proposition}
	\label{prop:small}
	Let $G$ be a graph.
	\begin{enumerate}
		\setlength{\itemsep}{0pt}
		\item $\gcg(G) = 1 \iff \Delta(G) = |V(G)| - 1$,
		\item $\gcg'(G) = 1 \iff$ $G$ is complete,
		\item $\gcg(G) = 2 \iff$ $G$ is the join of two non-complete graphs,
		\item $\gcg'(G) = 2 \iff$ $G$ is the join of two graphs, where at least one of them is non-complete.
	\end{enumerate}
\end{proposition}

\proof
Parts 1 and 2 are obvious. 

To prove part 3 notice that it holds that $\gcg(G) = 2$ if and only if $G$ does not have a universal vertex and there exists a vertex $u \in V(G)$ such that for every $v \in N(u)$ it holds $N[u] \cup N[v] = V(G)$. The latter is equivalent to the condition that $G$ is the join of two graphs which are not complete.

Similarly, $\gcg'(G) = 2$ if and only if $G$ is not complete and for every $u \in V(G)$ there exists its neighbor $v$ such that $N[u] \cup N[v] = V(G)$, which is equivalent to the fact that $G$ is the join of two graphs, one of them is non-complete.
\qed

Similarly as in~\cite{domSmall} for the domination game, we can find the relation between diameter of the graph and its connected game domination number.

\begin{proposition}
	\label{prop:diam}
	For a graph $G$ it holds that 
	$$\diam(G) \leq \gcg(G) + 1 \qquad \text{and} \qquad \diam(G) \leq \gcg'(G).$$
\end{proposition}

\proof
Let $\diam(G) = k$ and let $P$ be the longest induced path in $G$, $P = v_0 v_1 \ldots v_k$. Let $C$ be the set of vertices played in the course of a connected domination game on $G$. 

As $C$ is a dominating set of $G$, there exists a vertex $x \in C$ such that $x = v_0$ or $x \sim_G v_0$. Similarly, there exists a vertex $y \in C$ such that $y = v_0$ or $y \sim_G v_k$. As no path between $v_0$ and $v_k$ is shorter than $P$, it holds that $d_G(x,y) \geq k-2$. But we also have $d_C(x,y) \geq d_G(x,y)$. As $C$ is connected, it follows that $|C| \geq k-1$. Hence, $\gcg(G) \geq \diam(G) - 1$.

To prove the result for the S-game consider the following strategy of Staller. She starts the game on $v_0$ and then plays optimally. As the whole path $P$ must be dominated at the end of the game by a connected set, at least $k$ vertices are played. Hence, $\gcg'(G) \geq \diam(G)$.
\qed

Note that the bound in Proposition~\ref{prop:diam} is best possible, as $\gcg(P_n) = n-2 = \diam(P_n)-1$ and $\gcg'(P_n) = n-1 = \diam(P_n)$. On the other hand, the connected game domination number cannot be bounded from above by the diameter of a graph. The following proposition shows that given a diameter $r$ there exists a graph with arbitrary large connected game domination number. Note that if $n_1, \ldots, n_r \geq 2$, then $\diam(K_{n_1} \cp \cdots \cp K_{n_r}) = r$.

\begin{proposition}
	\label{prop:diamMax}
	Let $n_1, \ldots, n_r$ be integers, $2 \leq n_1 \leq \cdots \leq n_r$ and $n_r \geq 2 n_1 \cdots n_{r-1}$. It holds that
	$$\gcg(K_{n_1} \cp \cdots \cp K_{n_r}) = 2 n_1 \cdots n_{r-1}-1 \quad \text{and} \quad \gcg'(K_{n_1} \cp \cdots \cp K_{n_r}) = 2 n_1 \cdots n_{r-1}-2.$$
\end{proposition}

\proof
Iteratively using~\cite[Theorem 5]{connected} we have $\gcg(K_{n_1} \cp \cdots \cp K_{n_r}) \leq 2 n_1 \cdots n_{r-1}-1$. To prove the lower bound, consider the following strategy of Staller. She always plays in a copy of $K_{n_r}$ where Dominator made his first move. This is legal until the game is finished and assures that at least $\min\{n_r, 2 n_1 \cdots n_{r-1}-1\} = 2 n_1 \cdots n_{r-1}-1$ moves are played.

Similarly, we prove the result for the S-game. If Dominator plays all his moves in the copy of $K_{n_1} \cp \cdots \cp K_{n_{r-1}}$ where Staller made her first move, he assures that $\gcg'(K_{n_1} \cp \cdots \cp K_{n_r}) \leq 2 n_1 \cdots n_{r-1}-2$. If Staller plays her first move on an arbitrary vertex and all her other moves in the copy of $K_{n_r}$ where Dominator made his first move, she assures that $\gcg'(K_{n_1} \cp \cdots \cp K_{n_r}) \geq 2 n_1 \cdots n_{r-1}-2$.
\qed

\section{The upper bound for Staller-start connected game domination number}
\label{sec:staller-start}
%%%%%%%%%%%%%%%%%%%%%%%%%%%%%%%%%%%%%%%%%%%%%%%%%%%%%%%%%%
%%%%%%%%%%%%%%%%%%%%%%%%%%%%%%%%%%%%%%%%%%%%%%%%%%%%%%%%%%

For the classical domination game, we have $|\gg(G) - \gg'(G)| \leq 1$~\cite{dom, extremal}. An analogous result also holds for the total domination game~\cite{totDom}. For the connected domination game it holds that $$\gcg'(G) \geq \gcg(G) - 1$$ for all graphs $G$~\cite[Theorem 8]{connected}. Moreover, the bound is tight. The bound in the other direction, was not previously known. Below we present the solution to the following problem.

\begin{problem}[\cite{connected}, Problem 1]
	Find the maximum $k$ for which there exists a graph $G$ satisfying $\gcg'(G) = \gcg(G) + k$.
\end{problem}

In the following we show that such $k$ does not exist. However, for a fixed graph it can be bounded by $\gcg(G)$.

\begin{theorem}
	\label{thm:staller-start}
	For every graph $G$ it holds that $$\gcg'(G) \leq 2 \gcg(G).$$
\end{theorem}

\proof
Let $\gcg(G) = k$ and let $C$ be the set of played vertices in a D-game on $G$. Thus $C$ is connected, and for every vertex $v \in V(G)$ it holds that $v \in C$ or $v$ has a neighbor in $C$. Consider the following strategy for Dominator in the S-game on $G$ (the worst case scenario). Suppose Staller starts on $s_1 \in V(G) \setminus C$ and that the game is not over until all vertices in $C$ are played. Dominator's first move is in $C$ (on a neighbor of $s_1$). Even if Staller can always make a move outside of $C$, Dominator can play all the vertices there (due to connectedness) to finish the game at least after his $k$th move. Hence, at most $2k$ moves are played in total.
\qed

\begin{figure}[!ht]
	\begin{center}
		\begin{tikzpicture}[thick]
		
		\node[label=below: {$u_0$}] (u0) at (0,0) {};
		\node[label=135: {$u_1$}] (u1) at (1,0) {};
		\node[label=-135: {$u_2$}] (u2) at (2,0) {};
		\node[label=135: {$u_3$}] (u3) at (3,0) {};
		\node[label=-135: {$u_4$}] (u4) at (4,0) {};
		\node[label=135: {$u_5$}] (u5) at (5,0) {};
		\node[label=below: {$u_6$}] (u6) at (6,0) {};
		
		\node[label=above: {$x_1$}] (x1) at (1,1) {};
		\node[label=below: {$x_2$}] (x2) at (2,-1) {};
		\node[label=above: {$x_3$}] (x3) at (3,1) {};
		\node[label=below: {$x_4$}] (x4) at (4,-1) {};
		\node[label=above: {$x_5$}] (x5) at (5,1) {};
		
		\node[label=above: {$y_1$}] (y1) at (2,1) {};
		\node[label=below: {$y_2$}] (y2) at (3,-1) {};
		\node[label=above: {$y_3$}] (y3) at (4,1) {};
		\node[label=below: {$y_4$}] (y4) at (5,-1) {};
		\node[label=above: {$y_5$}] (y5) at (6,1) {};
		
		\node[label=above: {$z_1$}] (z1) at (2.5,1) {};
		\node[label=below: {$z_2$}] (z2) at (3.5,-1) {};
		\node[label=above: {$z_3$}] (z3) at (4.5,1) {};
		\node[label=below: {$z_4$}] (z4) at (5.5,-1) {};
		\node[label=above: {$z_5$}] (z5) at (6.5,1) {};
		
		\draw (u0) -- (u1) -- (u2) -- (u3) -- (u4) -- (u5) -- (u6);
		
		\draw (u1) -- (x1) -- (y1) -- (z1) -- (u2);
		\path (u2) edge (x1);
		\path (u2) edge (y1);
		\path (u2) edge (z1);
		
		\draw (u2) -- (x2) -- (y2) -- (z2) -- (u3);
		\path (u3) edge (x2);
		\path (u3) edge (y2);
		\path (u3) edge (z2);
		
		\draw (u3) -- (x3) -- (y3) -- (z3) -- (u4);
		\path (u4) edge (x3);
		\path (u4) edge (y3);
		\path (u4) edge (z3);
		
		\draw (u4) -- (x4) -- (y4) -- (z4) -- (u5);
		\path (u5) edge (x4);
		\path (u5) edge (y4);
		\path (u5) edge (z4);
		
		\draw (u5) -- (x5) -- (y5) -- (z5) -- (u6);
		\path (u6) edge (x5);
		\path (u6) edge (y5);
		\path (u6) edge (z5);

		\end{tikzpicture}
		\caption{The graph $G_6$.}
		\label{fig:primer}
	\end{center}
\end{figure}
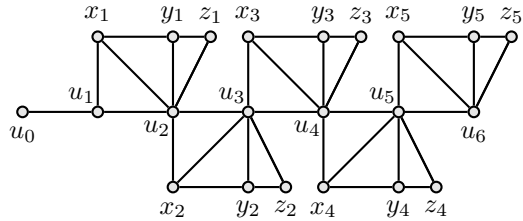

We now construct a graph that attains the equality in Theorem~\ref{thm:staller-start}. Let $G_n$, $n \geq 2$, be a graph with vertices $$V(G_n) = \{ u_0, \ldots, u_n, x_1, y_1, z_1, \ldots, x_{n-1}, y_{n-1}, z_{n-1} \}$$ and edges $u_i \sim u_{i+1}$ for $i \in \{0, \ldots, n-1\}$, $u_i \sim x_i \sim y_i \sim z_i$ and $u_{i+1} \sim x_i, y_i, z_i$ for $i \in \{1, \ldots, n-1\}$. See Figure~\ref{fig:primer} for the graph $G_6$.

\begin{lemma}
	\label{Gn-Dom}
	For $n \geq 2$ it holds $\gcg(G_n) = n$.
\end{lemma}

\proof
First notice that $\gc(G_n) = n$, as $\{u_1, \ldots, u_n\}$ is the smallest connected dominating set. Hence $\gcg(G_n) \geq \gc(G_n) = n$ by~\cite[Theorem 1]{connected}.

Observe the following strategy of Dominator. Dominator starts on $u_{n}$. The only playable vertex is $u_{n-1}$, hence Staller plays it. Now Dominator plays $u_{n-2}$ and Staller is forced to reply on $u_{n-3}$. This continues until the vertex $u_1$ is played and the game is finished after the $n$th move. Hence $\gcg(G_n) \leq n$.
\qed

\begin{lemma}
	\label{Gn-St}
	For $n \geq 2$ it holds $\gcg'(G_n) = 2n$.
\end{lemma}

\proof
It follows from Theorem~\ref{thm:staller-start} and Lemma~\ref{Gn-Dom} that $\gcg'(G_n) \leq 2 n$. Consider the following strategy of Staller to show that the equality is attained. But first set $I_1 = \{u_1, x_1\}$ and $I_j = \{y_{j-1}, u_j, x_j\}$ for $j \in \{2, \ldots, n-1\}$. 

Staller starts the game on $u_0$ and thus Dominator can only reply on $u_1$. From now on, Staller's strategy is to play on $V(G) \setminus \{u_2, \ldots, u_n\}$ whenever possible. By this she assures that at least two vertices are played on each $I_j$, $j \in \{1, \ldots, n-1\}$. 

Staller's second move is $x_1$, thus two moves are made on $I_1$. We now prove that at least two moves are made on $I_k$, $k \in \{2, \ldots, n-2\}$. If Dominator plays $u_k$ in the course of the game, then the vertices $x_k$, $y_k$, $z_k$, $u_{k+1}$, $x_{k+1}$, $y_{k+1}$, $z_{k+1}$, $u_{k+2}, \ldots$ have not been played yet (as Staller started on $u_0$ and the set of all played vertices must be connected). Hence Staller can reply on $x_k$, making a second move on $I_k$. 

If Staller plays $u_k$ in the course of the game, then the vertices $x_k$, $y_k$, $z_k$, $u_{k+1}$, $x_{k+1}$, $y_{k+1}$, $z_{k+1}$, $u_{k+2}, \ldots$ have not been played yet. But as Staller's strategy is to avoid playing on $\{u_2, \ldots, u_n\}$, this means that she was forced to play on $u_k$. Thus she was not able to play on $\{x_{k-1}, y_{k-1}, z_{k-1}\}$, which could only happen if Dominator played $y_{k-1}$ just before her move. Hence Staller's move on $u_k$ was already the second move on $I_k$. 

Additionally, to dominate the vertex $z_{n-1}$, at least one vertex from $\{u_n, y_{n-1}, z_{n-1}\}$ has to be played before the game ends. Hence, considering the first move of Staller on $u_0$, at least two moves on each $I_k$ and the last move to dominate $z_{n-1}$, this strategy assures that at least $1 + 2 (n-1) + 1 = 2n$ moves are played.
\qed

%%%%%%%%%%%%%%%%%%%%%%%%%%%%%%%%%%%%%%%%%%%%%%%%%%%%%%%%%%
%%%%%%%%%%%%%%%%%%%%%%%%%%%%%%%%%%%%%%%%%%%%%%%%%%%%%%%%%%
\section{Lexicographic products}
\label{sec:lexicographic}
%%%%%%%%%%%%%%%%%%%%%%%%%%%%%%%%%%%%%%%%%%%%%%%%%%%%%%%%%%
%%%%%%%%%%%%%%%%%%%%%%%%%%%%%%%%%%%%%%%%%%%%%%%%%%%%%%%%%%

In this section we present the connected game domination numbers in D- and S-games on the lexicographic product of graphs (for definition see Section~\ref{sec:intro}). Note that no such result is known for the (total) domination game. We first prove a lemma that will be useful in both games. We say that a first move played on some copy of $H$ is a \emph{new move} and all next moves on the same copy of $H$ are called \emph{duplicate moves}. Note that the first move on the graph $G[H]$ is always a new move, and if $\gcg(H) \geq 2$, then the second move in the game can be a duplicate move.

\begin{lemma}
	\label{lema:lexi}
	In a connected domination game on $G[H]$ it holds that after two new moves  are played, no duplicate moves are possible.
\end{lemma}

\proof
Let $S \subseteq V(G)$, $|S| \geq 2$, be a set of vertices in $G$ such that a new move was played on the copies $H_v$ for all $v \in S$. Clearly, $S$ is a connected set. Hence all vertices in copies $H_x$ for $x \in N[S]$ are already dominated, thus no other vertex in the copies $H_v$, $v \in S$, is playable.
\qed

To prove the general result for $G[H]$ we need to introduce another variation of the connected domination game -- a game when Staller skips her first move, i.e.\ Dominator plays two moves, and only then the players start to alternate moves. We denote the number of moves in such a game, when both players play optimally, with $\sgcg(G)$, and call the game the \emph{Staller-first-skip connected domination game}. Analogously, we define the \emph{Dominator-first-skip connected domination game} as the Staller-start game in which Dominator skips his first move, i.e., first Staller plays two moves and after that the players alternate taking moves. The number of moves in such a game when both players play optimally is $\sgcg'(G)$.

\begin{lemma}
	\label{lema:St-pass-1}
	For a graph $G$ it holds that \begin{enumerate}
		\setlength{\itemsep}{0pt}
		\item[(i)] $\gcg(G) - 1 \leq \sgcg(G) \leq \gcg(G) + 1,$
		\item[(ii)] $\gcg'(G) - 1 \leq \sgcg'(G) \leq \gcg'(G) + 1.$
	\end{enumerate}

\end{lemma}

\proof
We first prove (i). To prove the upper bound, consider the following strategy of Dominator. While the Staller-first-skip game is played on $G$, he imagines a connected domination game with Chooser on $G$. He plays his optimal first move in the imagined game and copies the move to the real game. In the real game he plays any legal move and copies it to the imagined game as Chooser's move. From now on he copies each move of Staller from the real to the imagined game, selects his optimal reply in the imagined game and copies it to the real game. Throughout this process the same set of vertices has been played in both games, thus all the copied moves are legal. As Staller played optimally in the real game and the same number of moves was played in the real and in the imagined game, it holds that $\sgcg(G)$ is at most the number of moves in a connected domination game on $G$ in which Chooser picks one vertex. Hence by the Chooser Lemma~\ref{lema:chooser}, $\sgcg(G) \leq \gcg(G) + 1$.

Similarly, we prove the lower bound with the following strategy of Staller. Parallel to the real Staller-first-skip game played on $G$, she imagines a D-game with Chooser on $G$. She copies the first move of Dominator from the real game to the imagined game as his first move there. But she copies his second move in the real game as a move of Chooser in the imagined game. From here on, she replies optimally in the imagined game and copies her move to the real game. She also copies each Dominator's move in the real game to the imagined game. Similarly as above, we use the Chooser Lemma and conclude that $\sgcg(G) \geq \gcg(G) - 1$.

Part (ii) is proved analogously.
\qed

The following examples (which have already been studied in~\cite{connected}) demonstrate that the bounds from Lemma~\ref{lema:St-pass-1}(i) cannot be improved.

\begin{example}
	\label{ex:st-pass}
	\begin{enumerate}
		\item For $n\geq3$ it clearly holds that $\gcg(P_n) = \sgcg(P_n) = n-2.$
		
		\item Let $F_1$ be the fan with $n \geq 7$ vertices and let $F_i$ be obtained from $F_{i-1}$ by identifying one of the vertices of degree $2$ and its neighbor of degree $3$ in $F_1$ and $F_{i-1}$, see Figure~\ref{fig:F_i}. Recall from~\cite{connected} that for $i \geq 1$ and all $n \geq 7$, it holds that $\gcg(F_i) = 2i-1$. Staller's strategy in a Staller-pass game in $F_i$ from~\cite{connected} also implies that $\sgcg(F_i) = 2i = \gcg(F_i) + 1$. 
		
		\begin{figure}[!hbt]
			\begin{center}
				\begin{tikzpicture}[thick,scale=1]
				\begin{scope}
				\node (a0) at (0,0) {};
				\node (a1) at (-1,0) {};
				\node (a2) at (-1,0.5) {};
				\node (a3) at (-0.6,0.8) {};
				\node (a4) at (0,1) {};
				\node (a5) at (0.6,0.8) {};
				\node (a6) at (1,0.5) {};
				\node (a7) at (1,0) {};
				
				\draw (a1) -- (a2) -- (a3) -- (a4) -- (a5) -- (a6) -- (a7);
				\path (a0) edge (a1);
				\path (a0) edge (a2);
				\path (a0) edge (a3);
				\path (a0) edge (a4);
				\path (a0) edge (a5);
				\path (a0) edge (a6);
				\path (a0) edge (a7);
				\end{scope}
				
				\begin{scope} [xshift=2.5cm]
				\node (a0) at (0,0) {};
				\node (a1) at (-1,0) {};
				\node (a2) at (-1,0.5) {};
				\node (a3) at (-0.6,0.8) {};
				\node (a4) at (0,1) {};
				\node (a5) at (0.6,0.8) {};
				\node (a6) at (1,0.5) {};
				\node (a7) at (1,0) {};
				
				\draw (a1) -- (a2) -- (a3) -- (a4) -- (a5) -- (a6) -- (a7);
				\path (a0) edge (a1);
				\path (a0) edge (a2);
				\path (a0) edge (a3);
				\path (a0) edge (a4);
				\path (a0) edge (a5);
				\path (a0) edge (a6);
				\path (a0) edge (a7);
				
				\node (b0) at (2,0) {};
				\node (b1) at (1,0) {};
				\node (b2) at (1,0.5) {};
				\node (b3) at (1.4,0.8) {};
				\node (b4) at (2,1) {};
				\node (b5) at (2.6,0.8) {};
				\node (b6) at (3,0.5) {};
				\node (b7) at (3,0) {};
				
				\draw (b1) -- (b2) -- (b3) -- (b4) -- (b5) -- (b6) -- (b7);
				\path (b0) edge (b1);
				\path (b0) edge (b2);
				\path (b0) edge (b3);
				\path (b0) edge (b4);
				\path (b0) edge (b5);
				\path (b0) edge (b6);
				\path (b0) edge (b7);
				\end{scope}
				
				\begin{scope} [xshift=7cm]
				\node (a0) at (0,0) {};
				\node (a1) at (-1,0) {};
				\node (a2) at (-1,0.5) {};
				\node (a3) at (-0.6,0.8) {};
				\node (a4) at (0,1) {};
				\node (a5) at (0.6,0.8) {};
				\node (a6) at (1,0.5) {};
				\node (a7) at (1,0) {};
				
				\draw (a1) -- (a2) -- (a3) -- (a4) -- (a5) -- (a6) -- (a7);
				\path (a0) edge (a1);
				\path (a0) edge (a2);
				\path (a0) edge (a3);
				\path (a0) edge (a4);
				\path (a0) edge (a5);
				\path (a0) edge (a6);
				\path (a0) edge (a7);
				
				\node (b0) at (2,0) {};
				\node (b1) at (1,0) {};
				\node (b2) at (1,0.5) {};
				\node (b3) at (1.4,0.8) {};
				\node (b4) at (2,1) {};
				\node (b5) at (2.6,0.8) {};
				\node (b6) at (3,0.5) {};
				\node (b7) at (3,0) {};
				
				\draw (b1) -- (b2) -- (b3) -- (b4) -- (b5) -- (b6) -- (b7);
				\path (b0) edge (b1);
				\path (b0) edge (b2);
				\path (b0) edge (b3);
				\path (b0) edge (b4);
				\path (b0) edge (b5);
				\path (b0) edge (b6);
				\path (b0) edge (b7);
				
				\node (c0) at (4,0) {};
				\node (c1) at (3,0) {};
				\node (c2) at (3,0.5) {};
				\node (c3) at (3.4,0.8) {};
				\node (c4) at (4,1) {};
				\node (c5) at (4.6,0.8) {};
				\node (c6) at (5,0.5) {};
				\node (c7) at (5,0) {};
				
				\draw (c1) -- (c2) -- (c3) -- (c4) -- (c5) -- (c6) -- (c7);
				\path (c0) edge (c1);
				\path (c0) edge (c2);
				\path (c0) edge (c3);
				\path (c0) edge (c4);
				\path (c0) edge (c5);
				\path (c0) edge (c6);
				\path (c0) edge (c7);
				\end{scope}

				\end{tikzpicture}
				\caption{Graphs $F_1$, $F_2$ and $F_3$.}
				\label{fig:F_i}
			\end{center}
		\end{figure}
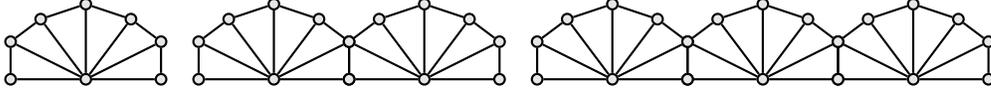
		
		\item Let $H_0$ be the graph from the left of Figure~\ref{fig:H_i}. Let $H_i$ be obtained by identifying two vertices in $H_0$ and $H_{i-1}$ in the same manner as above with $F_i$ (see Fig.~\ref{fig:H_i}). From~\cite{connected} it follows that for all $i \geq 1$, $\gcg(H_i) = 4i+2$ and from Chooser's strategy in~\cite{connected} we can conclude that $\sgcg(H_i) = 4i+1 = \gcg(H_i) - 1$. 
		
		\begin{figure}[!hbt]
			\begin{center}
				\begin{tikzpicture}[thick,scale=1]
				\begin{scope}
				\node (a0) at (0,0) {};
				\node (a1) at (-1,0) {};
				\node (a2) at (-1,0.5) {};
				\node (a3) at (-0.6,0.8) {};
				\node (a4) at (0,1) {};
				\node (a5) at (0.6,0.8) {};
				\node (a6) at (1,0.5) {};
				\node (a7) at (1,0) {};
				\node (ah) at (-0.4,1.2) {};
				
				\draw (a1) -- (a2) -- (a3) -- (a4) -- (a5) -- (a6) -- (a7);
				\path (a0) edge (a1);
				\path (a0) edge (a2);
				\path (a0) edge (a3);
				\path (a0) edge (a4);
				\path (a0) edge (a5);
				\path (a0) edge (a6);
				\path (a0) edge (a7);
				\path (ah) edge (a3);
				\path (ah) edge (a4);
				\end{scope}
				
				\begin{scope} [xshift=2.5cm]
				\node (a0) at (0,0) {};
				\node (a1) at (-1,0) {};
				\node (a2) at (-1,0.5) {};
				\node (a3) at (-0.6,0.8) {};
				\node (a4) at (0,1) {};
				\node (a5) at (0.6,0.8) {};
				\node (a6) at (1,0.5) {};
				\node (a7) at (1,0) {};
				\node (ah) at (-0.4,1.2) {};
				
				\draw (a1) -- (a2) -- (a3) -- (a4) -- (a5) -- (a6) -- (a7);
				\path (a0) edge (a1);
				\path (a0) edge (a2);
				\path (a0) edge (a3);
				\path (a0) edge (a4);
				\path (a0) edge (a5);
				\path (a0) edge (a6);
				\path (a0) edge (a7);
				\path (ah) edge (a3);
				\path (ah) edge (a4);
				
				\node (b0) at (2,0) {};
				\node (b1) at (1,0) {};
				\node (b2) at (1,0.5) {};
				\node (b3) at (1.4,0.8) {};
				\node (b4) at (2,1) {};
				\node (b5) at (2.6,0.8) {};
				\node (b6) at (3,0.5) {};
				\node (b7) at (3,0) {};
				\node (bh) at (1.6,1.2) {};
				
				\draw (b1) -- (b2) -- (b3) -- (b4) -- (b5) -- (b6) -- (b7);
				\path (b0) edge (b1);
				\path (b0) edge (b2);
				\path (b0) edge (b3);
				\path (b0) edge (b4);
				\path (b0) edge (b5);
				\path (b0) edge (b6);
				\path (b0) edge (b7);
				\path (bh) edge (b3);
				\path (bh) edge (b4);
				\end{scope}
				
				\begin{scope} [xshift=7cm]
				\node (a0) at (0,0) {};
				\node (a1) at (-1,0) {};
				\node (a2) at (-1,0.5) {};
				\node (a3) at (-0.6,0.8) {};
				\node (a4) at (0,1) {};
				\node (a5) at (0.6,0.8) {};
				\node (a6) at (1,0.5) {};
				\node (a7) at (1,0) {};
				\node (ah) at (-0.4,1.2) {};
				
				\draw (a1) -- (a2) -- (a3) -- (a4) -- (a5) -- (a6) -- (a7);
				\path (a0) edge (a1);
				\path (a0) edge (a2);
				\path (a0) edge (a3);
				\path (a0) edge (a4);
				\path (a0) edge (a5);
				\path (a0) edge (a6);
				\path (a0) edge (a7);
				\path (ah) edge (a3);
				\path (ah) edge (a4);
				
				\node (b0) at (2,0) {};
				\node (b1) at (1,0) {};
				\node (b2) at (1,0.5) {};
				\node (b3) at (1.4,0.8) {};
				\node (b4) at (2,1) {};
				\node (b5) at (2.6,0.8) {};
				\node (b6) at (3,0.5) {};
				\node (b7) at (3,0) {};
				\node (bh) at (1.6,1.2) {};
				
				\draw (b1) -- (b2) -- (b3) -- (b4) -- (b5) -- (b6) -- (b7);
				\path (b0) edge (b1);
				\path (b0) edge (b2);
				\path (b0) edge (b3);
				\path (b0) edge (b4);
				\path (b0) edge (b5);
				\path (b0) edge (b6);
				\path (b0) edge (b7);
				\path (bh) edge (b3);
				\path (bh) edge (b4);
				
				\node (c0) at (4,0) {};
				\node (c1) at (3,0) {};
				\node (c2) at (3,0.5) {};
				\node (c3) at (3.4,0.8) {};
				\node (c4) at (4,1) {};
				\node (c5) at (4.6,0.8) {};
				\node (c6) at (5,0.5) {};
				\node (c7) at (5,0) {};
				\node (ch) at (3.6,1.2) {};
				
				\draw (c1) -- (c2) -- (c3) -- (c4) -- (c5) -- (c6) -- (c7);
				\path (c0) edge (c1);
				\path (c0) edge (c2);
				\path (c0) edge (c3);
				\path (c0) edge (c4);
				\path (c0) edge (c5);
				\path (c0) edge (c6);
				\path (c0) edge (c7);
				\path (ch) edge (c3);
				\path (ch) edge (c4);
				\end{scope}

				\end{tikzpicture}
				\caption{Graphs $H_0$, $H_1$ and $H_2$.}
				\label{fig:H_i}
			\end{center}
		\end{figure}
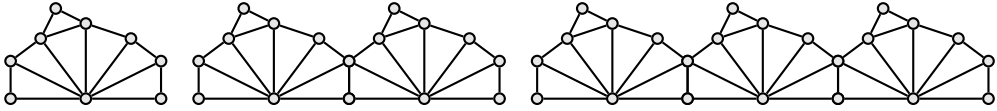
	\end{enumerate}
\end{example}

\begin{theorem}
	\label{thm:lexi-Dom}
	Let $G, H$ be graphs. It holds
	$$\gcg(G[H]) = \begin{cases}
	\sgcg(G) + 1; & \gcg(H) \geq 2 \text{ and } |V(G)| \geq 2,\\
	\gcg(G); & \gcg(H) = 1,\\
	\gcg(H); & |V(G)| = 1.
	\end{cases}$$
\end{theorem}

%PREVERI SE ENKRAT! (manjka $2; \gcg(G) = 1, |V(G)| \geq 2  \text{ and } \gcg(H) \geq 2$? Ne, je v 1.primeru.)

\proof
First consider the case $\gcg(H) \geq 2$ and $|V(G)| \geq 2$. Dominator's strategy is to start on a vertex in $H_{d_1}$, where $d_1$ is his optimal first move from the D-game on $G$. If Staller replies with a new move, then Lemma~\ref{lema:lexi} states that a game on $G[H]$ is essentially just a normal $D$-game on $G$. Additionally, by Lemma~\ref{lema:St-pass-1}(i) $\gcg(G) \leq \sgcg(G) + 1$. If Staller replies with a duplicate move on $H_{d_1}$, then Dominator's strategy is to reply on $H_{d_2}$, where $d_2$ is the optimal second move from the Staller-first-skip connected domination game on $G$. After this move, no more duplicate moves are possible. Hence at most $\sgcg(G)+1$ moves are played. 

On the other hand, Staller's strategy is to reply to the Dominator's first move with a duplicate move (she can as $\gcg(H) \geq 2$). If Dominator's second move is a new move, then they play at least $\sgcg(G) + 1$ moves. If Dominator's second move is a duplicate move, then Staller's strategy is to reply with a new move (with her optimal strategy from a D-game on $G$). Thus at least $\gcg(G) + 2 \geq \sgcg(G) + 1$ moves are played. This proves the first case.

Next, consider the case $\gcg(H) = 1$. As $\gcg(H) = 1$, after Dominator's first move on $h \in H_{d_1}$, where $d_1$ is Dominator's optimal first move on $G$ and $h$ is his optimal first move on $H$, Staller can only reply with a new move. And from here on, only new moves are possible due to Lemma~\ref{lema:lexi}. Hence, just a normal D-game is played, so $\gcg(G[H]) = \gcg(G)$.

Lastly, consider the case $|V(G)| = 1$. In this case $G[H] \cong H$, hence $\gcg(G[H]) = \gcg(H)$.
\qed

\begin{corollary}
	\label{cor:lexi-Dom}
	 Let $G$ and $H$ be graphs. If $\gcg(H) \geq 2$ and $|V(G)| \geq 2$, then
	 $$\gcg(G) \leq \gcg(G[H]) \leq \gcg(G) + 2.$$
\end{corollary}

Taking a graph $H$ with $\gcg(H) \geq 2$ and graphs from Example~\ref{ex:st-pass} as $G$ shows that all three possible values in Corollary~\ref{cor:lexi-Dom} can be achieved.

Now we focus on the S-game.

\begin{theorem}
	\label{thm:lexi-St}
	Let $G, H$ be graphs. It holds
	$$\gcg'(G[H]) = \begin{cases}
	\gcg'(G); & \gcg'(G) \geq 2,\\
	2; & \gcg'(G) = 1, |V(G)| \geq 2  \text{ and } \gcg'(H) \geq 2,\\
	\gcg'(H); & |V(G)| = 1, \text{ or } \gcg'(G) = 1, |V(G)| \geq 2  \text{ and } \gcg'(H) = 1.
	\end{cases}$$
\end{theorem}

%PREVERI SE ENKRAT!

\proof
First consider the case $\gcg'(G) \geq 2$. Dominator's strategy is to reply to the Staller's first move by playing a new vertex (he can as $\gcg'(G) \geq 2$). Lemma~\ref{lema:lexi} claims that in this case, just new moves are possible until the end of the game, thus a normal S-game is played on $G$. Hence, $\gcg'(G[H]) \leq \gcg'(G)$.  Staller's strategy is to start on $H_{s_1}$, where $s_1$ is her optimal first move on $G$. If Dominator plays a duplicate move (and if he can do so), she replies with a new move. Thus the number of moves is at least $\sgcg'(G) + 1 \geq \gcg'(G)$ by Lemma~\ref{lema:St-pass-1}(ii). If Dominator replies with a new move, the total number of moves is at least $\gcg'(G)$.

Next, consider the case $\gcg'(G) = 1$, $|V(G)| \geq 2$ and $\gcg'(H) \geq 2$. Dominator can reply on a new vertex and finish the game in the second move. Staller can start on $h \in H_{s_1}$ where $h$ is her optimal first move on $H$. Thus she forces at least two moves. 

Next, in the case $\gcg'(G) = 1$, $|V(G)| \geq 2$  and $\gcg'(H) = 1$, wherever Staller starts, she dominates the whole graph. But as $1 = \gcg'(H)$, this case can be joined with the case $|V(G)| = 1$ which is clear.
\qed

%%%%%%%%%%%%%%%%%%%%%%%%%%%%%%%%%%%%%%%%%%%%%%%%%%%%%%%%%%
%%%%%%%%%%%%%%%%%%%%%%%%%%%%%%%%%%%%%%%%%%%%%%%%%%%%%%%%%%
\section{Predomination}
\label{sec:predomination}
%%%%%%%%%%%%%%%%%%%%%%%%%%%%%%%%%%%%%%%%%%%%%%%%%%%%%%%%%%
%%%%%%%%%%%%%%%%%%%%%%%%%%%%%%%%%%%%%%%%%%%%%%%%%%%%%%%%%%

The behavior of domination games when one or more vertices are considered predominated has been widely studied, see for example~\cite{DomCritical, totDomCritical, ladders, totPredomination}. Let $S \subseteq V(G)$. When a domination game is played on the graph $G$ with vertices $S$ predominated, the vertices from $S$ can still be played (if the moves are legal), but they are already considered dominated. We denote such \emph{predominated graph} $G|S$. If $S = \{v\}$, we write $G|v$. A well known fact that holds for the classical domination game is the Continuation Principle~\cite{extremal}, which states the following. For a graph $G$ and sets $B \subseteq A \subseteq V(G)$ it holds that $\gg(G|A) \leq \gg(G|B)$ and $\gg'(G|A) \leq \gg'(G|B)$. Analogous results hold for the total domination game~\cite{totDom} as well as for other variations of the game~\cite{domGames}.

We first observe that a similar statement cannot hold for the connected domination game. For example, consider the path on $n \geq 5$ vertices and let $A$ be the set of all degree $2$ vertices on the path. As the played vertices must be connected, the connected domination game on $G|A$ can never be finished. 

However, it is still interesting to observe the case when only one vertex is predominated. In this case Dominator can achieve a finite number of moves in a D-game by starting the game on the predominated vertex. But as opposed to the other domination games, it seems that it might not be true that $\gcg(G|v) \leq \gcg(G)$. We return to this phenomenon in Proposition~\ref{prop:predom_increase}.

Note that in the S-game Staller might be able to prevent the game from ending even if just one vertex is predominated. For example, if a degree $2$ vertex on a path is predominated and Staller starts the game on a vertex at distance $2$ from it, the game can never end. Hence in the following, we only focus on the D-game on graphs with one vertex predominated.

\begin{proposition}
	\label{prop:predom_increase}
	There exists a graph $G$ with a vertex $v$ such that $$\gcg(G|v) = \gcg(G) + 1.$$
\end{proposition}

\proof
Consider the graph $G$ as presented on Figure~\ref{fig:G_1}. Due to the connectedness, the vertices $a,b,c,d,e,f,g$ must be played in any connected domination game on $G$. Thus $\gcg(G) \geq 7$. Consider the following strategy of Dominator to prove the equality. Let $d_1 = g$, then Staller can only reply on $s_1 = f$, Dominator strategy is to play $d_2 = e$, and then the remaining moves of the game are $s_2 = d$, $d_3 = c$, $s_3 = b$ and $d_4 = a$. Hence $\gcg(G) = 7$.

\begin{figure}[!hbt]
	\begin{center}
		\begin{tikzpicture}[thick,scale=1]
		\node[label=below: {$c$}, fill=black] (v) at (0,0) {};
		\node[label=below: {$b$}] (a1) at (-1,0) {};
		\node[label=below: {$a$}] (a2) at (-2,0) {};
		\node[label=below: {$a'$}] (a3) at (-3,0) {};
		\node[label=below: {$d$}] (b) at (1,0) {};
		\node[label=below: {$e$}] (c) at (2,0) {};
		\node[label=above: {$e'$}] (d) at (2,1) {};
		\node[label=below: {$f$}] (e) at (3,0) {};
		\node[label=above: {$f'$}] (f) at (3,1) {};
		\node[label=below: {$g$}] (g) at (4,0) {};
		\node[label=below: {$g'$}] (h) at (5,0) {};
		
		\draw (h) -- (g) -- (e) -- (c) -- (b) -- (v) -- (a1) -- (a2) -- (a3);
		\draw (c) -- (d) -- (f) -- (e);

		\end{tikzpicture}
		\caption{Graph $G$.}
		\label{fig:G_1}
	\end{center}
\end{figure}

Consider the D-game played on $G|c$. Dominator starts the game on $N[c]$ (otherwise the game does not end). From a simple case analysis it follows that in all these three cases Staller will be able to make a move on $e'$, hence $\gcg(G|c) = 8 = \gcg(G) + 1$.
\qed

Next, consider some simple examples. Let $n \geq 3$ and $P_n$ be a path with vertices $1,\ldots, n$ and naturally defined edges. We know that $\gcg(P_n) = n-2$~\cite{connected}. If $v \in \{1,n\}$, then it clearly holds $\gcg(P_n|v) = n-3$. However, if $v \in \{2,\ldots,n-1\}$, then we have $\gcg(P_n|v) = n-2$ (as the set of played vertices must be connected). This reasoning can be extended to trees~\cite{connected}. For a tree $T$ it holds that $\gcg(T) = \gc(T) = \gcg(T|v)$ for every vertex $v \in V(T)$ that is not a leaf. 

Let $n \geq 4$ and $C_n$ be a cycle with vertices $1,\ldots, n$ and naturally defined edges. Clearly, it holds that $\gcg(C_n) = n-2$. Let $v \in V(C_n)$ and consider the D-game on $C_n|v$. Without loss of generality $v = 1$. Wherever Dominator starts the game, at least the vertices $3, \ldots, n-1$ or any set of $n-2$ vertices must be played during the game. Thus $\gcg(C_n|v) \geq n-3$. Dominator's strategy is to start on $3$. From this move on, the moves in the game are uniquely determined ($s_1 = 4$, $d_2 = 5$ \ldots). Hence, $\gcg(C_n|v) = n-3$.

Note that on the contrary to the previous example of cycles where for every vertex $v \in V(G)$ it holds that $\gcg(G|v) < \gcg(G)$, it cannot happen that $\gcg(G|v) > \gcg(G)$ for all vertices $v \in V(G)$. Let $d_1$ be an optimal first move of Dominator in a D-game on $G$. Then $\gcg(G|d_1) \leq \gcg(G)$ as Dominator can still start the game on $d_1$.

But there are also examples where predominating a vertex cannot decrease the connected domination game number.

\begin{proposition}
	\label{prop:cut-vertex}
	Let $G$ be a graph with a cut vertex $u$ (i.e.\ $G - u$ is not connected). It holds that $\gcg(G|u) \geq \gcg(G)$.
\end{proposition}

\proof
Let $\gcg(G) = k$. Suppose $\gcg(G|u) \leq k-1$. Let $S$ be the set of played vertices in the strategy that achieves this bound. If the vertices from $S$ dominate $u$, then this same strategy can be used on $G$ to show that $\gcg(G) \leq k-1$, which is not the case. Thus the vertices from $S$ do not dominate $u$, meaning that $S \cap N[u] = \emptyset$. But then $S$ cannot be connected, a contradiction.
\qed

Two other families of graphs with the same property as cycles are circular ladders and M\"{o}bius ladders. Predomination in the total domination game on those families has already been studied in~\cite{ladders}. Recall that the \emph{circular ladder} $\CL_n$, $n \geq 3$, is the Cartesian product $C_n \cp K_2$. Denote the vertices of $CL_n$ by $\{(i, j) \; ; \; i \in [n], j \in [2] \}$ and by $C_n^1$, $C_n^2$ the two $C_n$-layers in $\CL_n$. The \emph{M\"{o}bius ladder} $\ML_n$ is a graph obtained from $C_{2n}$ by adding edges between opposite vertices.

\begin{theorem}
	\label{thm:circular}
	For $n \geq 4$ and every vertex $v \in V(\CL_n)$ it holds that $$\gcg(\CL_n) = 2 (n-2) \quad \text{and} \quad \gcg(\CL_n|v) = 2 (n-2) - 1.$$
\end{theorem}

\proof First we prove that $\gcg(\CL_n) = 2 (n-2) = 2n - 4$. To prove the upper bound, consider the following strategy of Dominator. His first move is $(1,1)$ and as long as possible he plays on $V(C_n^1) \setminus \{ (n-1,1), (n,1), (1,1) \}$. If Staller makes no move on $C_n^1$, her first move can only be on $(1,2)$, so she dominates the vertex $(n,2)$. After Dominator's $(n-2)$-th move at most one vertex can be undominated, i.e.\ $(n-1,2)$. So in this case Staller finishes the game and at most $2 (n-2)$ moves are played. If Staller makes at least one move on $C_n^1$, then after the players play $n-2$ moves on $C_n^1$ two additional moves might be needed (to dominate $(n-1,2)$ and $(n,2)$). Thus the number of moves is at most $2(n-3) + 2 = 2n - 4$. 

%To prove the upper bound suppose that after $2 n - 4$ moves were played, at least one vertex remains undominated. It is easy to see that the only configuration with four unplayed vertices and at least one vertex not yet dominated is (without loss of generality) $(n,1)$, $(1,1)$, $(2,1)$, $(1,2)$. Suppose vertex $(2,2)$ was played earlier in the game that the vertex $(n,2)$. As the set of played vertices must be connected and all other vertices were played, $(n,2)$ can only be played last. But after all vertices except $(n,2)$, $(n,1)$, $(1,1)$, $(2,1)$ and $(1,2)$ were played, $(n,2)$ is not a legal move. As we reached a contradiction it follows that $\gcg(\CL_n) \leq 2n - 4$. 

To prove the lower bound consider the following strategy of Staller. When Dominator plays on $(i, j)$, Staller replies on $(i, 3-j)$. Clearly, such reply is always legal and the game finishes after at least $2 (n-2)$ moves are played. Thus $\gcg(\CL_n) \geq 2n-4$.

Next we prove that $\gcg(\CL_n|v) = 2 (n-2)-1 = 2n - 5$. Without loss of generality let $v = (1,1)$. 

To prove the upper bound consider the following strategy of Dominator. His first move is $(3,2)$ and as long as possible he plays on $V(C_n^2) \setminus \{ (1,2), (2,2), (3,2) \}$. If Staller makes no move on $C_n^2$, her first move can only be on $(3,1)$, so she dominates the vertex $(2,1)$. The game is over after Dominator's $(n-2)$-th move (when he plays $(n,2)$). Thus the number of moves is at most $(n-2) + (n-3) = 2n - 5$. If Staller makes at least one move on $C_n^2$, then after the players play $n-2$ moves on $C_n^2$ one additional move might be needed (to dominate $(2,1)$). Thus the number of moves is at most $2(n-3) + 1 = 2n - 5$. 

To prove the lower bound consider the following strategy of Staller. When Dominator plays on $(i, j)$, Staller replies on $(i, 3-j)$. Clearly, such reply is legal unless the game is over, so the game finishes after at least $2 (n-2) - 1$ moves are played. Thus $\gcg(\CL_n) \geq 2n-5$.
\qed

Notice that the strategies described in the proof of Theorem~\ref{thm:circular} are also valid on M\"{o}bius ladders.

\begin{corollary}
	\label{cor:mobius}
	For $n \geq 4$ and every vertex $v \in V(\ML_n)$ it holds that $$\gcg(\ML_n) = 2 (n-2) \quad \text{and} \quad \gcg(\ML_n|v) = 2 (n-2) - 1.$$
\end{corollary}

%%%%%%%%%%%%%%%%%%%%%%%%%%%%%%%%%%%%%%%%%%%%%%%%%%%%%%%%%%
%%%%%%%%%%%%%%%%%%%%%%%%%%%%%%%%%%%%%%%%%%%%%%%%%%%%%%%%%%
\section{Concluding remarks}
\label{sec:concluding}
%%%%%%%%%%%%%%%%%%%%%%%%%%%%%%%%%%%%%%%%%%%%%%%%%%%%%%%%%%
%%%%%%%%%%%%%%%%%%%%%%%%%%%%%%%%%%%%%%%%%%%%%%%%%%%%%%%%%%

Recall that for the domination game (resp.\ total domination game), critical graphs are defined as graphs $G$ such that for every $v \in V(G)$ it holds that $\gg(G|v) < \gg(G)$ (resp.\ $\gtg(G|v) < \gtg(G)$)~\cite{DomCritical, totDomCritical}. On the other hand, it is not clear how to define connected domination game critical graphs. From Section~\ref{sec:predomination} we know that there exist graphs with the property that $\gcg(G|v) < \gcg(G)$ for all vertices $v$, and that it can also happen that $\gcg(G|v) > \gcg(G)$, but this cannot be true on all vertices $v$. However, the following is not clear.

\begin{problem}
	\label{prob:critical}
	 Does there exist a graph $G$ with the following properties:
	 \begin{itemize}
	 	\setlength{\itemsep}{0pt}
	 	\item for all $v \in V(G)$: $\gcg(G|v) \neq \gcg(G)$, and
	 	\item there exists a vertex $u \in V(G)$: $\gcg(G|u) > \gcg(G)$.
	 \end{itemize}
\end{problem}

If the answer to the Problem~\ref{prob:critical} is positive, than a sensible definition for critical graphs could be the following. A graph $G$ is connected domination game critical if for every vertex $v \in V(G)$ it holds that $\gcg(G|v) \neq \gcg(G)$. However, if the answer to the Problem~\ref{prob:critical} is negative, a definition should go along the same lines as the definition of critical graphs for domination and total domination games. So a graph $G$ would be connected domination game critical if for every vertex $v \in V(G)$ it would hold that $\gcg(G|v) < \gcg(G)$.

Two other interesting problems arise from the predomination of vertices.

\begin{problem}
	\label{prob:pmk}
	\begin{enumerate}
		\setlength{\itemsep}{0pt}
		\item What is the maximal value of $k$ such that there exists a graph $G$ with a vertex $v \in V(G)$ with the property $\gcg(G|v) = \gcg(G) + k$?
		\item What is the maximal value of $k$ such that there exists a graph $G$ with a vertex $v \in V(G)$ with the property $\gcg(G|v) = \gcg(G) - k$?
	\end{enumerate}
\end{problem}

Recall that for domination and total domination game the answers to analogous problems is known and equals $k = 0$ for the first part and $k = 2$ for the second part~\cite{DomCritical, totPredomination}.

\vspace{1.5cc}
\noindent \textbf{Acknowledgements.}
The author would like to thank Sandi Klav\v zar for many helpful conversations and suggestions.


\begin{thebibliography}{99}
	
\bibitem{connected}
M.\ Borowiecki, A.\ Fiedorowicz, E.\ Sidorowicz,
Connected domination game,
Applicable Analysis and Discrete Mathematics, to appear,
Doi: https://doi.org/10.2298/AADM171126020B.

\bibitem{domGames}
B.\ Bre\v{s}ar et al.,
The variety of domination games, 
arXiv:1807.02695 [math.CO].
	
%\bibitem{EdgeVertexRemoval}
%B.\ Bre\v{s}ar, P.\ Dorbec, S.\ Klav\v{z}ar, G.\ Ko\v{s}mrlj, Domination game: effect of edge- and vertex-removal, Discrete Math.\ 330 (2014) 1--10.

%\bibitem{pspace}
%B.\ Bre\v{s}ar, M.\ A.\ Henning, The game total domination problem is log-complete in PSPACE, Inform.\ Process.\ Lett.\ 126 (2017) 12--17.

\bibitem{dom}
B.\ Bre\v{s}ar, S.\ Klav\v{z}ar, D.\ F.\ Rall, Domination game and an imagination strategy, SIAM J.\ Discrete Math.\ 24 (2010) 979--991.

\bibitem{3/4_2}
Cs.\ Bujt\'{a}s, M.\ A.\ Henning, Z.\ Tuza, Transversal game on hypergraphs and the $\frac{3}{4}$-Conjecture on the total domination game, SIAM J.\ Discrete Math.\ 30 (2016) 1830--1847.

\bibitem{DomCritical}
Cs.\ Bujt\'{a}s, S.\ Klav\v zar, G.\ Ko\v smrlj, Domination game critical graphs, Discuss.\ Math.\ Graph Theory 35 (2015) 781--796.

\bibitem{totDomCP}
P.\ Dorbec, M.\ A.\ Henning, Game total domination for cycles and paths, Discrete Appl.\ Math.\ 208 (2016) 7--18.

%\bibitem{union}
%P.\ Dorbec, G.\ Ko\v{s}mrlj, G.\ Renault, The domination game played on unions of graphs, Discrete Math.\ 338 (2015) 71--79.

\bibitem{G-e}
M.\ A.\ Henning, W.\ B.\ Kinnersley, Domination game: A proof of the $3/5$-conjecture for graphs with minimum degree at least two, SIAM J.\ Discrete Math.\ 30 (2016) 20--35.

\bibitem{totDom}
M.\ A.\ Henning, S.\ Klav\v zar, D.\ F.\ Rall, Total version of the domination game, Graphs Combin.\ 31 (2015) 1453--1462.

\bibitem{totDomCritical}
M.\ A.\ Henning, S.\ Klav\v zar, D.\ F.\ Rall, Game total domination critical graphs, Discrete Appl.\ Math.\ 250 (2018) 28--37.

\bibitem{3/4_1}
M.\ A.\ Henning, S.\ Klav\v zar, D.\ F.\ Rall, The $4/5$ upper bound on the game total domination number, Combinatorica 37 (2017) 223--251.

\bibitem{ladders}
M.\ A.\ Henning, S.\ Klav\v zar, Infinite families of circular and M\"{o}bius ladders that are total domination game critical, Bull.\ Malays.\ Math.\ Sci.\ Soc.\ 41 (2018) 2141--2149

\bibitem{dom35}
M.\ Henning, C.\ L\"{o}wenstein,
Domination game: extremal families for the 3/5-conjecture for forests,
Discuss.\ Math.\ Graph Theory 37 (2017) 369--381.

\bibitem{3/4_3}
M.\ A.\ Henning, D.\ F.\ Rall, Progress towards the total domination game $\frac{3}{4}$-Conjecture, Discrete Math.\ 339 (2016) 2620--2627.

\bibitem{trees}
M.\ A.\ Henning, D.\ F.\ Rall, Trees with equal total domination and game total domination numbers, Discrete Appl.\ Math.\ 226 (2017) 58--70.

\bibitem{totPredomination}
V.\ Ir\v{s}i\v{c}, Effect of predomination and vertex removal on the game total domination number of a graph, Discrete Appl.\ Math.\ (2018), https://doi.org/10.1016/j.dam.2018.09.011.

\bibitem{book}
T.\ James, P.\ Dorbec, A.\ Vijayakumar (2017) Further Progress on the Heredity of the Game Domination Number. In: S.\ Arumugam, J.\ Bagga, L.\ Beineke, B.\ Panda (eds) Theoretical Computer Science and Discrete Mathematics. ICTCSDM 2016. Lecture Notes in Computer Science, vol 10398. Springer, Cham.

\bibitem{extremal}
W.\ B.\ Kinnersley, D.\ B.\ West, R.\ Zamani, Extremal problems for game domination number, SIAM J.\ Discrete Math.\ 27 (2013) 2090--2017.

\bibitem{domSmall}
S.\ Klav\v{z}ar, G.\ Ko\v{s}mrlj, S.\ Schmidt,
On graphs with small game domination number,
Appl.\ Anal.\ Discrete Math.\ 10 (2016) 30--45.

\bibitem{domEdgeCuts}
S.\ Klav\v{z}ar, D.\ F.\ Rall,
Domination game and minimal edge cuts,
Discrete Math.\ 342 (2019) 951--958.

\bibitem{domPathsCycles}
G.\ Ko\v{s}mrlj,
Domination game on paths and cycles,
Ars Math.\ Contemp.\ 13 (2017) 125--136.
	
%
%\bibitem{D-trivial}
%M.\ J.\ Nadjafi-Arani, M.\ Siggers, H.\ Soltani, Characterisation of forests with trivial game domination numbers, J.\ Comb.\ Optim.\ 32 (2016), no.\ 3, 800--811.

\bibitem{domLargest}
K.\ Xu, X.\ Li, S.\ Klav\v{z}ar,
On graphs with largest possible game domination number,
Discrete Math.\ 341 (2018) 1768--1777.
\end{thebibliography}
\end{document}